\magnification=1200
\font\bb=msbm10

\def\sq{\hbox{\rlap{$\sqcap$}$\sqcup$}}

\centerline{\bf PATH-COMPONENT INVARIANTS FOR SPACES } 
\centerline{\bf OF POSITIVE SCALAR CURVATURE METRICS}
\bigskip
\centerline{\bf David J. Wraith}
\bigskip\medskip
\leftskip 1cm \rightskip 1cm \noindent {\bf Abstract:} {\it The Kreck-Stolz $s$-invariant is a classic path-component invariant for the space of positive scalar curvature metrics on certain spin manifolds, with $|s|$ an invariant of the path-component in the moduli space. It is an absolute (as opposed to relative) invariant, but this strength comes at the expense of being defined only under restrictive topological conditions. The aim of this paper is to construct an analogous invariant for certain product manifolds on which the $s$-invariant is not defined.}
\bigskip\medskip
\centerline{\bf \S0 Introduction}
\bigskip
\leftskip 0cm \rightskip 0cm 
Given a manifold $M$ which supports a positive scalar curvature metric, an important but difficult question is to determine what can be said about the topology of the space of all positive scalar curvature metrics on $M$, $\hbox{Riem}_{scal \ge 0}(M).$ (In this paper we will always assume that spaces of metrics are equipped with the smooth topology.) One can also ask about the topology of the {\it moduli} space of positive scalar curvature metrics on $M$. Recall that the diffeomorphism group $\hbox{Diff}(M)$ acts on the space of all metrics $\hbox{Riem}(M)$ by pull-back, and this action preserves the property of positive scalar curvature. Thus we can form the moduli space of positive scalar curvature metrics $$\hbox{Riem}_{scal \ge 0}(M)/\hbox{Diff}(M).$$ One can pose analogous questions for other curvature conditions, such as positive or non-negative Ricci curvature, negative sectional curvature etc. There has been much recent activity in this general direction: for example see [BHSW], [BERW], [HSS], [CS], [CM], [Wa1], [Wa2], [Wr1], [Wr2], [BH], [DKT], [FO1-3] and the book [TW]. 

In this paper we will focus on spaces of positive scalar curvature metrics. However we wish to highlight out at the outset that using the results of [Wr2], all statements involving positive scalar curvature can be easily modified to yield analogous statements about non-negative scalar curvature. Although these corresponding results are stronger, we have chosen to focus on positivity in order to simplify the exposition. 

A basic tool for studying spaces of positive scalar curvature metrics is the Kreck-Stolz $s$-invariant (see [KS] or [TW] for details). Under the appropriate topological conditions this allows one to distinguish between different path components of the space of positive scalar curvature metrics, and even between path-components of the moduli space of positive scalar curvature metrics. As we will need to refer to these conditions regularly, for convenience we make the following definition:
\proclaim Definition 0.1. A closed spin manifold $M$ of dimension $4k-1$, $k \ge 2$, which admits a positive scalar curvature metric and for which all real Pontrjagin classes vanish will be said to satisfy the {\it Kreck-Stolz conditions}.
\par
It was shown in [KS] that if $M$ satisfies the Kreck-Stolz conditions and $(M,g)$ has positive scalar curvature, then $s(M,g)$ is an invariant of the path-component of of positive scalar curvature metrics containing $g$. Moreover, if $H^1(M;\hbox{\bb Z}_2)=0$ (which ensures that given an orientation for $M$ the spin structure is unique), then $|s|$ can be shown to be an invariant of the path-component of the {\it moduli space} of positive scalar curvature metrics containing $[g]$. This was used in [KS] to show that the moduli space of positive scalar curvature metrics for any manifold $M$ with $H^1(M;\hbox{\bb Z}_2)=0$ satisfying the Kreck-Stolz conditions has infinitely many path-components. (In stark contrast, it was shown in [CM] that the moduli space of positive scalar curvature metrics on closed orientable 3-manifolds is path-connected, provided this space is non-empty.) Using the obvious fact that metrics of positive sectional or positive Ricci curvature also have positive scalar curvature, by examining the underlying space of positive scalar curvature metrics using the $s$-invariant Kreck and Stolz were able to show that in dimension seven there are manifolds for which the moduli space of positive Ricci curvature metrics has infinitely many path-components, and also examples with positive sectional curvature for which the moduli space of such metrics is not path-connected. The $s$-invariant has subsequently been used to establish analogous results in other contexts. For example the author showed in [Wr1] that the moduli space of Ricci positive metrics on all homotopy spheres in dimensions $4n-1 \ge 7$ which bound a parallelisable manifold has infinitely many path-components, showing that this infinite disconnectedness phenomenon occurs through an infinite range of dimensions, and providing the first examples away from dimension seven. Very recently, in [DKT] it is shown that in every dimension $4n-1 \ge 7$, there are infinitely many closed manifolds for which the moduli space of non-negative sectional curvature metrics has infinitely many path-components.

To provide some context, it has long been known (see [LM; IV Theorem 7.7]) that for any closed spin manifold $M$ of dimension $4n-1\ge 7$ admitting positive scalar curvature, $\hbox{Riem}_{scal>0}(M)$ has infinitely many path-components. It was pointed out in [PS; Remark 2.26] that the same argument used to establish [LM; IV Theorem 7.7] can also be used to show that the corresponding moduli space has infinitely many path-components. In the author's experience this second point is not so widely known. In both cases, the argument can be expressed neatly using Gromov and Lawson's relative index. This is defined for pairs of positive scalar curvature metrics $g_0,g_1$ on $M$, and is given by $$i(g_0,g_1)=\hbox{ind}D^+(M \times [0,1],g),$$ where $D^+$ denotes the Dirac operator and $g$ is any metric on $M \times [0,1]$ restricting to $dt^2+g_0$ and $dt^2+g_1$ in a neighbourhood of the boundary components. It can be shown that this is an invariant of the path-components of positive scalar curvature metrics to which $g_0$ and $g_1$ belong, and vanishes if both belong to the same component. The advantage of the Kreck-Stolz $s$-invariant over this is that it is an {\it absolute} invariant, i.e. it only depends on a single metric. Indeed $i(g_0,g_1)=s(M,g_0)-s(M,g_1)$ whenever the right-hand side is defined. However Kreck and Stolz show ([KS; 2.16]) that it is not possible to define an absolute invariant of this type without imposing extra topological conditions on $M$.

The aim of this paper is to demonstrate that it is possible to make similar constructions under alternative topological circumstances to those in Definition 0.1. We achieve this by providing an extension of the $s$-invariant to certain product manifolds. The new setting is as follows: we consider Riemannian product manifolds $(M,g_M) \times (N,g_N),$ where $M$ satisfies the Kreck-Stolz conditions, $g_M$ has positive scalar curvature, and $N$ is a closed spin manifold of dimension $4l$, $l \ge 1$, with $\hat{A}(N) \neq 0.$

For manifolds in dimensions congruent to 0 modulo 4 the $\hat{A}$-genus is a topological obstruction to the existence of positive scalar curvature metrics. Nevertheless, any Riemannian product involving a positive scalar curvature metric on one factor can be adjusted by scaling to produce a positive scalar curvature metric. In the above product, there is some very small $c>0$ such that the metric $c^2g_M+g_N$ has positive scalar curvature.

The key point here is that the $s$-invariant is not defined for product manifolds of this type. To see this consider $\hat{A}(N).$ The $\hat{A}$-genus is a rational linear combination of rational Pontrjagin numbers, and hence if $\hat{A}(N) \neq 0,$ this means that some real Pontrjagin class of $N$ is non-zero, and in turn this means that some real Pontrjagin class of $M \times N$ is also non-zero. Thus the Kreck-Stolz conditions are not satisfied by the product $M \times N$ in this case. 

Let us summarise our new context in a definition:
\proclaim Definition 0.2. A closed oriented Riemannian spin manifold $(X,g)$ with positive scalar curvature will be said to have a Kreck-Stolz product structure if it is 
orientation preserving 
isometric to a Riemannian product manifold $(M^{4(k-l)-1},g_M) \times (N^{4l},g_N),$ $k-l\ge 2$, $l \ge 1$, where $M$ satisfies the Kreck-Stolz conditions (Definition 0.1), $\hat{A}(N) \neq 0,$ and $H^1(X;\hbox{\bb Z}_2)=0.$ In this case, the Kreck-Stolz product structure is a 7-tuple $(X,g,\phi,M,N,g_M,g_N),$ where $\phi:(X,g) \to (M \times N,g_M+g_N)$ is the orientation preserving isometry. We will denote the set of all Kreck-Stolz product structures on $X$ by ${\cal K}(X).$
\par
\noindent{\it Remark 1:} It is implicit in the above definition that the orientations on $M$ and $N$ are chosen so as to be compatible with $X$ under the isometry $\phi$. Thus a different choice of $\phi$ might result in the orientations on $M$ or $N$ having to be rechosen. Now it follows from the K\"unneth Theorem that the condition $H^1(X;\hbox{\bb Z}_2)=0$ forces both $H^1(M;\hbox{\bb Z}_2)=0$ and $H^1(N;\hbox{\bb Z}_2)=0,$ and thus a unique spin structure on $X$ means that the spin structures on $M$ and $N$ are also unique for a given orientation. Consequently, the fact that $\phi$ is orientation preserving automatically means that it is spin-structure preserving.
\smallskip
\noindent{\it Remark 2:} One can also make certain uniqueness statements about the factors $M$ and $N$ appearing in a Kreck-Stolz product structure: if $(X,g)$ admits a Kreck-Stolz product structure, then this structure is uniquely determined up to isometry of the individual product factors. This is discussed in detail in Corollary 2.16.
\smallskip
Notice that the diffeomorphism group of $X$ acts on ${\cal K}(X)$ is a natural way: for $\theta \in \hbox{Diff}(X)$ we set $$\theta \cdot (X,g,\phi,M,N,g_M,g_N)=(X,\theta^*(g),\phi\circ\theta,M,N,g_M,g_N).$$ Thus we can also consider the moduli space of such structures, ${\cal K}(X)/\hbox{Diff}(X).$

The main result we will establish in this paper as follows: 
\proclaim Theorem 0.3. Consider a Riemannian manifold $(X,g)$ which has a Kreck-Stolz product structure $(X,g,\phi,M,N,g_M,g_N).$ Then there is a function $\tilde{s}:{\cal K}(X)\to \hbox{\bb Q}$ such that given any other positive scalar curvature metric $g'$ on $X$ for which $(X,g')$ has a Kreck-Stolz product structure $(X,g',\phi,M,N,g'_M,g'_N),$ (so the diffeomorphism $\phi:X \to M \times N$ is shared by both structures), the following statements hold.
\item{(i)} The relative index $i(g,g')=\tilde{s}(X,g)-\tilde{s}(X,g'),$ so in particular if $g$ and $g'$ belong to the same path-component of the space of positive scalar curvature metrics on $X,$ then $\tilde{s}(X,g)=\tilde{s}(X,g').$
\item{(ii)} $\hbox{Riem}_{scal \ge 0}(X)$ has infinitely many path-components of positive scalar curvature metrics distinguished by $\tilde{s}.$
\item{(iii)} If $[g]$ denotes the class of $g$ in the moduli space of positive scalar curvature metrics on $X,$ then for any $h \in [g],$ $|\tilde{s}(X,h)|=|\tilde{s}(X,g)|$.
\item{(iv)} $|\tilde{s}|$ descends to give a {\bb Q}-valued function on the moduli space of Kreck-Stolz structures ${\cal K}(X)/\hbox{Diff}(X).$ 
\par
\noindent {\it Notational remark:} It might appear that $\tilde{s}$ depends on the whole Kreck-Stolz product structure, however this is not the case. By Remark 2 above, the Kreck-Stolz structure only depends on the pair $(X,g)$ up to orientation preserving isometry of the factors $(M,g_M)$ and $(N,g_N),$ and we will see in due course that $\tilde{s}$ is invariant under such isometries. Hence our notation $\tilde{s}(X,g)$ is justified.

Note that if we allowed the degenerate case $l=0$, i.e. the case where $N$ is a point, $\tilde{s}$ would reduce to $s$. 

To illustrate Theorem 0.3, we will present some explicit examples. These will take the form of product manifolds, for which we can take the diffeomorphism in the Kreck-Stolz structure to be the identity map.

Recall that a K3 surface $K^4$ satisfies $\hat{A}(K^4)=-2$, and so this simply-connected spin manifold cannot support a metric of positive scalar curvature. Similarly there is a simply-connected spin `Bott manifold' $B^8$ for which $\hat{A}(B^8)=1,$ so this too does not admit a positive scalar curvature metric. (The Bott manifold can be constructed by forming the boundary connected sum of 28 copies of the manifold constructed by plumbing the tangent disk bundle of $S^4$ to itself according to the $E_8$-graph. The resulting object has boundary $S^7,$ and this can then be made into a smooth closed manifold $B^8$ by gluing in a disc $D^8$. Together with $\hbox{\bb H}P^2,$ $B^8$ generates $\Omega_8^{spin}\cong \hbox{\bb Z} \oplus \hbox{\bb Z}.$) We note that despite the fact that neither $K^4$ nor $B^8$ admit positive scalar curvature, both are known to admit Ricci flat metrics. Using Theorem 0.3 together with the definition of $\tilde{s}$ (Definition 2.8) and results from [Wr1; page 2014] we immediately obtain:
\proclaim Theorem 0.4. If $K^4$ denotes the K3 surface, $B^8$ the Bott manifold, and $\Sigma^{4n-1}$ is any homotopy $n$-sphere ($n \ge 2$) which bounds a parallelisable manifold, then there is a sequence $g_j$ of Ricci positive metrics on $\Sigma$ such that given Ricci flat metrics $g_K$ on $K^4$ and $g_B$ on $B^8$ we have $$\eqalign{\tilde{s}(\Sigma \times K^4, g_j+g_K)&=-{{j|bP_{4n}|+q} \over {2^{2n-3}(2^{2n-1}-1)}}; \cr \tilde{s}(\Sigma \times B^8, g_j+g_B)&={{j|bP_{4n}|+q} \over {2^{2n-2}(2^{2n-1}-1)}}, \cr}$$ where $q=q(\Sigma)$ is an integer depending on $\Sigma$, and where $bP_{4n}$ denotes the group of diffeomorphism classes of homotopy spheres bounding a parallelisable manifold of dimension $4n$. In particular, for different $j$ the metrics $g_j+g_K$ respectively $g_j+g_B$ belong to different path-components of the space of positive scalar curvature metrics on $\Sigma \times K^4$ respectively $\Sigma \times B^8.$
\par
Theorem 0.4 should be compared with Theorem 0.7 in [Wr2]. In fact, combining Theorem 0.3 with the results in [Wr2] it is not difficult to show that the condition of positive scalar curvature in Theorem 0.4 can be replaced with non-negative Ricci curvature.

We should also point out that products involving the Bott manifold appear in the stable Gromov-Lawson-Rosenberg conjecture (see for example [S2;\S1], or [RS;4.17]). Specifically, this claims that the {\it Rosenberg index} of a connected closed spin manifold $M$ of dimension at least five vanishes if and only if for some $k \ge 1,$ the manifold $M \times (B^8)^k$ admits a positive scalar curvature metric. Here $(B^8)^k$ denotes the $k$-fold product.

Whereas the modulus of the Kreck-Stolz $s$-invariant gives a well-defined invariant on the path-components of the moduli space of positive scalar curvature metrics, we are not in a position to assert that the same is true for $\tilde{s}.$ The essential difference is that metrics supporting Kreck-Stolz product structures will typically be distributed discretely through a path component of positive scalar curvature metrics, and we can only make a direct comparison between those which share the same spin-structure preserving isometry to a product, up to a product isometry of the individual factors. It could be that two such metrics lie in different path-components of the space of positive scalar curvature metrics and yield different $|\tilde{s}|$ values, but some isometric copy of one of the metrics lies in the same path-component as the other. In this situation, both metrics give rise to equivalance classes in the moduli space of positive scalar curvature metrics which belong to the same path component, though their $|\tilde{s}|$-values disagree. This just serves to highlight the difficulties associated with trying to extend path-component invariants beyond the classical results.

Theorem 0.3 raises the obvious question as to whether $\tilde{s}$ is a path-component invariant for at least {\it some} Kreck-Stolz product manifolds $M \times N.$ It turns out that this question is intimately related to the question of which manifolds are $\hat{A}$-multiplicative fibres, in the sense of [HSS; Definition 1.8]. Of particular relevance here is whether or not such products can be $\hat{A}$-multiplicative fibres {\it in degree 0}, which means that for any fibre bundle $F \to E \to S^1$ where the fibre $F=M\times N,$ we have $\hat{A}(E)=0.$ (In [HSS; Proposition 1.9] it is shown that the vanishing of all rational Pontrjagin classes is sufficient to satisfy this condition, however this does not help for Kreck-Stolz products.) If $F$ is such a fibre, then for any diffeomorphism $\phi:F \to F$ and any positive scalar curvature metric $g$ on $F$, the Gromov-Lawson relative index $i(g,\phi^*g)=0.$ This is a consequence of the the fact that the Atiyah-Patodi-Singer index formula is well-behaved under gluing manifolds: if $(V_1,g_1),(V_2,g_2)$ are Riemannian spin manifolds with boundary, and $\psi:\partial V_1 \to -\partial V_2$ is a spin-structure preserving isometry, then $\hbox{ind}D^+_{V_1}+\hbox{ind}D^+_{V_2}=\hat{A}(V_1\cup_\psi V_2).$ The argument is as follows. Consider a metric on $F \times [0,1]$ interpolating between $g$ on one boundary component and $\phi^*g$ on the other. Assume as usual that the interpolating metric is a product near each boundary. Now form the mapping torus $T_\phi$ of the diffeomorphism $\phi:F \to F,$ noticing that our metric on $F \times [0,1]$ descends to give a well-defined metric $\bar{g}$ on $T_\phi.$ If we remove a small piece of the form $F \times [0,\epsilon]$ from this torus, where we are assuming that $\epsilon$ is so small that the metric on $F \times [0,1]$ is a product throughout $F \times [0,\epsilon],$ then clearly the relative index for the boundary metrics is $i(g,g)=0.$ The relative index of the remaining part of the mapping torus is $i(g,\phi^*g),$ but by the above gluing formula, $$\eqalign{i(g,\phi^*g)&=i(g,g)+i(g,\phi^*g)\cr &=\hbox{ind}D^+(F\times [0,\epsilon],dt^2+g)+\hbox{ind}D^+(T_\phi \setminus (F\times [0,\epsilon]),\bar{g})\cr &=\hat{A}(T_\phi).\cr}$$ Now if $F$ is an $\hat{A}$-multiplicative fibre in degree 0, then $\hat{A}(T_\phi)=0,$ and the relative index claim follows. (Alternatively, we could argue from [HSS; Proposition 2.2].)

The significance of the equation $i(g,\phi^*g)=0$ for $\tilde{s}$ is as follows. Suppose that $g,g'$ are positive scalar curvature metrics on a manifold $X$ which belong to the same path-component of positive scalar curvature metrics, and for which both $(X,g)$ and $(X,g')$ have Kreck-Stolz product structures with respect to isometries $\psi,\theta:X \to M \times N$ respectively. Thus $(X,\psi^*(\theta^{-1})^*g')$ has a Kreck-Stolz structure with isometry $\psi.$ By Theorem 0.3(i) we deduce that $$i(g,\psi^*(\theta^{-1})^*g')=\tilde{s}(X,g)-\tilde{s}(X,\psi^*(\theta^{-1})^*g'),$$ but this relative index vanishes if $X$ is an $\hat{A}$-multiplicative fibre. By Theorem 0.3(iii) we also have that $\tilde{s}(X,g')=\tilde{s}(X,\psi^*(\theta^{-1})^*g'),$ and consequently $\tilde{s}(X,g)=\tilde{s}(X,g').$ Allowing for the fact that diffeomorphisms might reverse orientation and change the sign of $\tilde{s}$ (see Lemma 2.11), in the case of an $\hat{A}$-multiplicative fibre, we deduce that $|\tilde{s}|$ is a path-component invariant for the moduli space of positive scalar curvature metrics on $X.$ 

Thus we are motivated to pose the following problem, in part to find examples of manifolds on which $|\tilde{s}|$ is a path-component invariant for the moduli space of positive scalar curvature metrics, and in part to extend the scope of [HSS; Proposition 1.9]:
\proclaim Problem 0.5. Find examples of Kreck-Stolz product manifolds $M \times N$ which are $\hat{A}$-multiplicative fibres in degree 0.
\par
In relation to the above problem, it is worth remarking that either such examples are commonplace, which would be interesting from the point-of-view of $\tilde{s},$ or they are not. In the latter case, this means that there must be rich classes of bundles in which the $\hat{A}$-genus does not behave multiplicatively. To the best of the author's knowledge, only one family of bundles in which $\hat{A}$ is not multiplicative is known at present, namely that constructed in [HSS; Theorem 1.4]. 

As $\tilde{s}$ is only defined for metrics isometric to a product, this naturally leads to the following question, which we believe is of independent interest:
\proclaim Question 0.6. Consider a product manifold $M \times N$ which admits a positive scalar curvature metric. Can one find conditions on $M$ and $N$ under which every path-component of positive scalar curvature metrics contains a product metric?
\par

Note that there are situations in which $M \times N$ admits positive scalar curvature metrics but no product metric with positive scalar curvature. For example consider the case where $M$ is a simply-connected spin 4-manifold with $\hat{A}(M)=0$ which does not admit a positive scalar curvature metric (see [R; Counterexample 1.13]), and where $N$ is a K3 surface. As noted above, a K3 surface is a simply-connected spin manifold with non-zero $\hat{A}$-genus, and therefore does not support a metric of positive scalar curvature. The product $M \times N$ is then a simply-connected spin 8-manifold with $\hat{A}(M \times N)=\hat{A}(M)\hat{A}(N)=0.$ By Gromov-Lawson [GL], all simply-connected spin 8-manifolds with vanishing $\hat{A}$-genus admit positive scalar curvature metrics. The author is grateful to Boris Botvinnik for pointing out this example. Of course, dimension four manifolds are somewhat special from a positive scalar curvature point of view. For an example involving higher-dimensional factors, one could start with the the manifold $M^5$ described by Thomas Schick in [S1], which is a counter-example to the (unstable) Gromov-Lawson-Rosenberg conjecture. (See for example [RS; \S4] for a general discussion on this.) Now $M^5$ does not admit a positive scalar curvature metric, however for some $k \ge 1$ the product $M^5 \times (B^8)^k$ {\it does} admit such a metric (where $B^8$ is the Bott manifold as above), as the stable Gromov-Lawson-Rosenberg conjecture is known to hold for $M^5$. 
Of course $M^5 \times (B^8)^n$ cannot admit a product metric with positive scalar curvature as none of the factors individually support such a metric.

In a slightly different direction we note that product metrics play a crucial role in the recent paper [TWi], which investiagtes the moduli space of non-negative Ricci curvature metrics on certain manifolds, and will similarly be a central feature in a forthcoming paper of Boris Botvinnik and the author on the same topic. 

This paper is laid out as follows. In \S1 we outline the construction of the Kreck-Stolz $s$-invariant, as this provides the blueprint for establishing Theorem 0.3. The construction of $\tilde{s}$ and the proof of Theorem 0.3 are contained in \S2. 

The author would like to thank Bernd Ammann for reading a preliminary version of this paper and for providing some useful comments. He is grateful to Anand Dessai and Wilderich Tuschmann for alerting him to [PS; Remark 2.26]. He would also like to thank Boris Botvinnik for useful conversations. Special thanks go to the anonymous referee who identified a significant issue with the previous version of this paper.

\bigskip\bigskip


\centerline{\bf \S1 The Kreck-Stolz $s$-invariant}
\bigskip

We begin by recalling the index theorem of Atiyah-Patodi-Singer: 
\proclaim Theorem 1.1. ([APS]) Let $(W,g_W)$ be a compact even dimensional Riemannian spin manifold with non-empty boundary $M$, where the metric $g_W$ is a product $dt^2+g_M$ in a neighbourhood of the boundary. Consider the Atiyah-Singer Dirac operator $D^+$ on $W$ acting on the subspace of spinor bundle sections for which the restriction to $M$ belongs to the span of the negative eigenspaces of the operator induced on $M$. Then the index of this (restricted domain) Dirac operator on $W$ is given by $$\hbox{ind}D^+(W,g_W)=\int_W \hat{A}(p_*(W,g_W)) -{{h(M,g_M)+\eta(M,g_M)} \over 2},$$ where $\hat{A}$ denotes the $\hat{A}$-polynomial in the Pontrjagin forms of the metric, $h$ is the dimension of the space of harmonic spinors on the boundary $M$, and $\eta$ is the eta-invariant of the Dirac operator on $M$.
\par

Later on we will encounter eta-invariants of other operators. We will adopt the convention that if an operator is not specified in the notation it should be assumed to be a Dirac operator.

The integral term appearing in Theorem 1.1 appears to depend both on the topology of $W$ and the metric on $W$. However it is not difficult to see that only the metric in a neighbourhood of the boundary actually influences the value of the integral. (The argument behind this is detailed in Lemma 2.5.) In the light of this observation it is natural to ask: can we separate out the topological dependence on $W$ from the metric dependence near the boundary? The answer to this is a qualified yes: the integral of any summand in the integrand can be rewritten in this way provided it is decomposable (in the sense that it involves a product of forms), and provided that $M$ has vanishing real Pontrjagin classes.

From now on let us assume that $M$ does indeed have vanishing real Pontrjagin classes. Let $\alpha,\beta$ denote Pontrjagin forms or products of Pontrjagin forms on $W$. As a consequence of the above assumption together with the product structure of the metric $g_W$ near the boundary, following the notation in [KS; 2.8] we can define a form $d^{-1}(\alpha \wedge \beta)$ on $M$ by setting $d^{-1}(\alpha \wedge \beta)=\hat{\alpha}\wedge (\beta|_M),$ where $\hat{\alpha}$ satisfies $d\hat{\alpha}=\alpha|_M.$ A simple Stokes' Theorem argument then shows that $$\int_W \alpha \wedge \beta =\int_M d^{-1}(\alpha \wedge \beta)\, +\, \langle j^{-1}[\alpha] \cup j^{-1}[\beta],[W,M]\rangle,$$ where $j:H^*(W,M;\hbox{\bb R}) \to H^*(W;\hbox{\bb R})$ is the map induced by inclusion, and the angled brackets denote evaluation on the fundamental homology class. By the cohomology long exact sequence of the pair $(W,M)$ it is easy to see that we need $M$ to have vanishing real Pontrjagin classes in order for the required pre-images under $j$ to be defined.

Suppose now that $W$ has dimension $4k$. The top-dimensional term of the $\hat{A}$-polynomial has all its summands decomposable except for the term in $p_k(W,g_W).$ In order to deal with this, a linear combination of the $\hat{A}$ and $L$-polynomials is formed which has zero $p_k$ term. Specifically, Theorem 1.1 is applied to $\hat{A}+a_kL$ where $a_k=1/(2^{2k+1}(2^{2k-1}-1)).$ The resulting index formula is given by $$\eqalign{\hbox{ind}D^+(W,g_W)=&\int_M d^{-1}(\hat{A}+a_kL)(p_*(M,g_M)) \cr -& {{h(M,g_M)+\eta(M,g_M)} \over 2} \cr -&a_k\eta(B(M,g_M))-t(W), \cr}$$ where $\eta(B(M,g_M))$ is the eta-invariant of the signature operator $B$ on $M$, and $t(W)$ is the topological term $$t(W)=-\langle (\hat{A}+a_kL)(j^{-1}p_*(W)),[W,M]\rangle+a_k\sigma(W)$$ where the $p_\ast(W)$ are the Pontrjagin classes of $W$ (as opposed to forms), and $\sigma(W)$ is the signature. If we assume that the scalar curvature of $M$ is positive, this forces $h(M,g_M)=0.$

The idea behind the $s$-invariant is to collect together all the terms depending on the boundary $(M,g_M)$:
\proclaim Definition 1.3. Given a closed spin manifold $M^{4k-1}$ with positive scalar curvature and vanishing real Pontrjagin classes, the $s$-invariant is given by $$s(M,g)=-{1 \over 2}\eta(M,g_M)-a_k\eta(B(M,g_M))+\int_M d^{-1}(\hat{A}+a_kL)(p_*(M,g_M)).$$
\par
With this definition we see immediately that $$s(M,g_M)=\hbox{ind}D^+(W,g_W)+t(W).$$ Moreover if the metric $g_W$ also has positive scalar curvature then the index term above vanishes, leaving $s(M,g_M)=t(W).$ In this situation $s$ is completely determined by the topology of the bounding manifold $W$.

Note that Definition 1.3 does not require $M$ to be the boundary of a suitable manifold $W$. 

The key properties of the $s$-invariant can be proved by applying the above analysis to the case $W=M \times I$ for any interval $I$, after showing that $s$ is additive across disjoint unions and is sign-sensitive to orientation. Specifically it can be shown (as already mentioned in \S0) that $s$ is a path-component invariant for the space of positive scalar curvature metrics. Furthermore if $H^1(M;\hbox{\bb Z}_2)=0$ then $|s(M,g)| \in \hbox{\bb Q}$ is an invariant of the path-component of the moduli space of positive scalar curvature metrics on $M$ containing $g$. (See [KS; Proposition 2.13].)  

\bigskip\bigskip


\centerline{\bf \S2 An extension of the Kreck-Stolz $s$-invariant}
\bigskip

Let $W^{4(k-l)}$ ($k>l \ge 1$) and $N^{4l}$ be compact oriented spin manifolds. We suppose that $W$ has boundary $M$ (possibly disconnected), and that $N$ is a closed manifold. If $\pi_W$ (respectively $\pi_N$) denote the projections of $W \times N$ onto $W$ (respectively $N$), then by the Whitney formula the total rational or real Pontrjagin class satisfies $p(W \times N)=p(\pi^*_W(TW))p(\pi^*_N(TN))$ since $T(W \times N) \cong \pi^*_W(TW) \oplus \pi^*_N(TN).$ By the multiplicative property of the $\hat{A}$-polynomial we obtain the equality of polynomials $\hat{A}(p(W \times N))=\hat{A}(p(\pi^*_W(TW)))\hat{A}(p(\pi^*_N(TN))).$ By the naturality of the Pontrjagin classes we can then write $$\hat{A}(p(W \times N))=\pi^*_W(\hat{A}(p(M))\pi^*_N(\hat{A}(p(N)).$$

\proclaim Lemma 2.1. With $W$, $N$ as above, suppose we choose a product metric $g_W+g_N$ on $W \times N$. Then for the Pontrjagin forms corresponding to this metric we have $$p_i(W \times N;g_W+g_N)=\sum_{j+k=i}\pi_W^* p_j(W;g_W)\wedge\pi_N^* p_k(N;g_N).$$
\par
\medskip
\noindent{\bf Proof.} The Pontrjagin forms are symmetric polynomials in the curvature form for the given metric. Recall that given a local tangent frame field $\{s_i\}$ for a Riemannian $n$-manifold, the curvature form $\Omega$ is an $(n \times n)$-matrix of 2-forms $(\Omega^i_j)$ with entries defined by $$R(X,Y)(s_j)=\sum_{i=1}^n \Omega^i_j(X,Y)s_i.$$ For a product metric such as $g_W+g_N$ and frame fields $s_1,...,s_{4(k-l)} \in \Gamma(TW \oplus 0) \subset T(W \times N)$ and $s_{4(k-l)+1},...,s_{4k} \in \Gamma(0\oplus TN),$ the curvature 2-form satisfies $$\Omega=\pmatrix{\Omega_W&0 \cr 0&\Omega_N}$$ where $\Omega_W$ and $\Omega_N$ are the pull-backs of the curvature forms of $(W,g_W)$ respectively $(N,g_N).$ (See [Mo] page 208.) The total Pontrjagin form is then given by $$\det \Bigl(\hbox{\bb I}-{1 \over {2\pi i}}\Omega\Bigr)=\det \Bigl(\hbox{\bb I}-{1 \over {2\pi i}}\Omega_W\Bigr)\wedge \det \Bigl(\hbox{\bb I}-{1 \over {2\pi i}}\Omega_N\Bigr),$$ where the determinants on the right-hand side are the pull-backs (to $W \times N$) of the total Pontrjagin forms of $(W,g_W)$ and $(N,g_N).$ The lemma then follows by expanding these total classes into their individual terms. \hfill \sq 
\medskip

\proclaim Corollary 2.2. With the set-up of Lemma 2.1 we have the following decomposition of $\hat{A}$-polynomials into Pontrjagin forms: $$\hat{A}(p_*(W \times N;g_W+g_N))=\pi^*_W\hat{A}(p_*(W;g_W)) \wedge \pi^*_N\hat{A}(p_*(N;g_N)).$$ In particular for the top-dimensional forms we have $$\hat{A}_k(p_*(W \times N;g_W+g_N))=\pi^*_W\hat{A}_{k-l}(p_*(W;g_W)) \wedge \pi^*_N\hat{A}_l(p_*(N;g_N)).$$
\par
\medskip
\noindent{\bf Proof.} As discussed above, the equivalent formula to the first statement holds for Pontrjagin classes, and follows from the decomposition of those classes on product manifolds. By Lemma 2.1 Pontrjagin forms for product manifolds equipped with product metrics decompose into terms involving the individual factor manifolds in exactly the same way as Pontrjagin classes. The result follows immediately. For the second statement we simply note that a top dimensional form on $W \times N$ can only be formed from a product of top dimensional forms on the factors, since any higher degree form must be zero. \hfill\sq

\proclaim Lemma 2.3. Given top-dimensional differential forms $\alpha$, $\beta$ on oriented manifolds $X$ respectively $Y$, we have $$\int_{X \times Y} \pi^*_X \alpha \wedge \pi^*_Y \beta=\Bigl(\int_X \alpha \Bigr)\Bigl(\int_Y \beta\Bigr).$$
\par
\medskip
\noindent{\bf Proof.} This equation holds as it holds locally in any coordinate neighbourhood which is a product of coordinate neighbourhoods for $X$ and $Y$ individually. Using such a coordinate system the calculation reduces to showing that for appropriate functions $a$ and $b$: $$\int_{U \times V} a(x_1,...,x_r)b(y_1,...,y_s) \,dx_1...dx_rdy_1...dy_s$$ $$=\int_U a(x_1,...,x_r)\,dx_1...dx_r \int_V b(y_1,...,y_s) \,dy_1...dy_s,$$ which holds trivially. \hfill \sq
\medskip

\proclaim Corollary 2.4. With $W$, $N$ and metrics as above we have $$\int_{W \times N} \hat{A}_k(p_*(W \times N;g_W+g_N))=\hat{A}(N)\int_W \hat{A}_{k-l}(p_*(W;g_W)),$$ where $\hat{A}(N)$ is the $\hat{A}$-genus of $N$.
\par
\medskip
\noindent{\bf Proof.} This follows immediately from Corollary 2.2, Lemma 2.3 and the fact that $\hat{A}(N)=\int_N \hat{A}_l(p_*(N;g_N)).$
\medskip

{\it From now on we will assume that the product metric $g_W+g_N$ takes the form $dt^2+g_M+g_N$ near the boundary.} Since the Pontrjagin forms are defined using the Levi-Civita connection of the metric, if we replace the $g_W+g_N$ by another metric $g$ which takes the same product form $dt^2+g_M+g_N$ near the boundary we will change the Pontrjagin forms, however this will not change the value of the integral over $W \times N$:

\proclaim Lemma 2.5. Consider an oriented manifold $X^{4n}$ with non-empty connected boundary $Y$. Let $\phi$ be a top dimensional Pontrjagin form (respectively a top dimensional wedge product of Pontrjagin forms) on $X$ corresponding to a Riemannian metric $g_X,$ which is a product $dt^2+g_Y$ near the boundary. If $\phi'$ is the top dimensional Pontrjagin form (respectively the corresponding wedge product of Pontrjagin forms) arising from a metric $g'_X$ which also takes the form $dt^2+g_Y$ near the boundary, then $$\int_X \phi=\int_X \phi'.$$ 
\par
\medskip
\noindent{\bf Proof.} We consider the metric $g_X \cup g'_X$ on the oriented double of $X$, $X \cup (-X).$ This metric is smooth as the individual metrics agree near the common boundary. Now all oriented double manifolds are oriented boundaries, and hence all Pontrjagin numbers of $X\cup (-X)$ must vanish. Thus if we let $\psi$ be the top dimensional Pontrjagin form (respectively wedge product of Pontrjagin forms) on $X\cup(-X)$ arising from $g_X \cup g'_X$, then $$\int_{X\cup(-X)} \psi =0.$$ But $$\eqalign{\int_{X\cup(-X)} \psi &=\int_X \phi\,+\,\int_{-X} \phi' \cr &=\int_X \phi\,-\,\int_X \phi'. \cr}$$ Thus $\int_X \phi=\int_X \phi'$ as claimed. \hfill \sq
\medskip
From Corollary 2.4 and Lemma 2.5 we obtain:
\proclaim Corollary 2.6. Let $W$, $N$, $g_W$ and $g_N$ be as before, with $g_W$ taking the form $dt^2+g_M$ near $\partial W=M,$ and let $g$ be any metric on $W \times N$ which takes the same product form $dt^2+g_M+g_N$ as $g_W+g_N$ near the boundary. Then $$\int_{W \times N} \hat{A}_k(p_*(W \times N;g))=\hat{A}(N)\int_W \hat{A}_{k-l}(p_*(W;g_W)),$$ and applying the Atiyah-Patodi-Singer index theorem to $(W \times N;g)$ we obtain $$\hbox{ind}D^+(W \times N;g)=\hat{A}(N)\int_W\hat{A}_{k-l}(p_*(W;g_W)) \, -{{h+\eta} \over 2}(M \times N;g_M+g_N).$$
\par
Following [KS], {\it from now on we will make the assumption that the real Pontrjagin classes of $M =\partial W$ vanish}. This assumption allows us to re-write the above integral. Following the argument and notation in [KS] as outlined in \S1 we obtain
\proclaim Proposition 2.7. With $W$, $N$ and $g$ as above, and assuming the real Pontrjagin classes of $M$ vanish, $$\eqalign{\hbox{ind}D^+(W \times N;g)=&\hat{A}(N)\Bigl[\int_M d^{-1}(\hat{A}+a_{k-l}L)(p_*(M;g_M)) \, -a_{k-l}\eta(B(M,g_M))-t(W)\Bigr] \cr &-{{h+\eta} \over 2}(M \times N;g_M+g_N), \cr}$$ where $L$ is the Hirzebruch $L$-polynomial, $B$ denotes the signature operator, $$a_n:=1/(2^{2n+1}(2^{2n-1}-1)),$$ and the topological term $t(W)$ is given by $$t(W)=-\Big\langle(\hat{A}+a_{k-l}L)(j^{-1}p_*(W)),[W,M]\Big\rangle+a_{k-l}\sigma(W)$$ where $j$ denotes the inclusion map $j:H^*(W,M;\hbox{\bb R}) \to H^*(W;\hbox{\bb R})$ and $\sigma(W)$ is the signature of $W$.
\par
\noindent{\bf Proof.} It follows from Theorem 1.1 and the Atiyah-Patodi-Singer index theorem applied to the signature operator ([APS; 4.14]) that $$\int_W (\hat{A}+a_{k-l}L)(p_*(W;g_W)) =\int_W \hat{A}(p_*(W,g_W)) + a_{k-l}\sigma(W)+a_{k-l}\eta(B(M,g_M)).$$ By [KS; Lemma 2.7] we have $$\eqalign{\int_W (\hat{A}+a_{k-l}L)(p_*(W;g_W)) =& \int_M d^{-1}(\hat{A}+a_{k-l}L)(p_*(M;g_M))\cr &+\Big\langle(\hat{A}+a_{k-l}L)(j^{-1}p_*(W)),[W,M]\Big\rangle. \cr}$$ Combining the above two statements with Corollary 2.6 yields the result. \hfill\sq

If we assume that both $g_M$ and $g_M+g_N$ have positive scalar curvature, (we can always achieve this by scaling $g_M$ if necessary), then the term $h(M \times N;g_M+g_N)$ in the statement of Proposition 2.7 is zero. Collecting together the boundary terms as in [KS] then leads to
\proclaim Definition 2.8(a). Given $M$ and $N$ as above (so in particular the real Pontrjagin classes of $M$ all vanish and $\hat{A}(N) \neq 0$), together with a positive scalar curvature metric $g_M$ on $M$ and a metric $g_N$ on $N$ such that the product metric $g_M +g_N$ on $M \times N$ has positive scalar curvature, we set $$\eqalign{\tilde{s}(M \times N,g_M+g_N):=&\hat{A}(N)\Bigl[\int_M d^{-1}(\hat{A}+a_{k-l}L)(p_*(M;g_M))-a_{k-l}\eta(B(M,g_M))\Bigr]\cr &-{1 \over 2}\eta(D_{(M\times N;g_M+g_N)}). \cr }$$
\par
\medskip
We now can write $$\hbox{ind}D^+(W \times N;g)=\tilde{s}(M \times N,g_M+g_N)-\hat{A}(N){t}(W). \eqno{(\dag)}$$ Recalling the definition of the $s$-invariant (Definition 1.3) allows us to re-express $\tilde{s}$ as $$\tilde{s}(M\times N,g_M+g_N)=\hat{A}(N)s(M,g_M)+{1 \over 2}\hat{A}(N)\eta(M,g_M)-{1 \over 2}\eta(M \times N;g_M+g_N).$$

If the metric $g$ on $W\times N$ has positive scalar curvature then the index term vanishes and we are left with
\proclaim Lemma 2.9. With all manifolds and metrics as above, if $g$ is a positive scalar curvature metric on $W \times N$ (which as always is a product $dt^2+g_M +g_N$ near the boundary) then $$\tilde{s}(M \times N,g_M+g_N)=\hat{A}(N)t(W).$$ The right-hand side of this expression depends only on the topology of $W\times N$, and is independent of the choice of metrics.
\par

We now point out some key properties of $\tilde{s}$. The arguments needed here are essentially the same as those required to establish the equivalent properties for $s$. Although these arguments are for the most part suppressed in [KS], they are explained in depth in Chapter 5 of [TW], and we therefore omit the details here. (In relation to Lemma 2.11 below, we remark that the oriented manifold $-(M \times N)$ can be viewed as either $(-M)\times N$ or $M \times (-N),$ with the same conclusion obtained in either case following the arguments on pages 45-46 of [TW].)
\proclaim Lemma 2.10. $\tilde{s}$ is additive over disjoint unions in the following sense: $$\tilde{s}((M_1 \times N)\sqcup (M_2 \times N),g_{M_1}+g_N \sqcup g_{M_2}+g_N)=\tilde{s}(M_1 \times N,g_{M_1}+g_N)+\tilde{s}(M_2 \times N,g_{M_2}+g_N).$$
\par

\proclaim Lemma 2.11. $\tilde{s}$ is sensitive to the orientation of $M$ in the sense that $\tilde{s}(M \times N,g_M+g_N)=-\tilde{s}(-(M \times N),g_M+g_N).$
\par

\proclaim Lemma 2.12. $\tilde{s}$ is additive over connected sums in the following sense: $$\tilde{s}((M_1 \sharp M_2) \times N),(g_{M_1}\sharp g_{M_2})+g_N)=\tilde{s}(M_1 \times N,g_{M_1}+g_N)+\tilde{s}(M_2 \times N,g_{M_2}+g_N),$$ where $g_{M_1} \sharp g_{M_2}$ is the (canonical) Gromov-Lawson positive scalar curvature metric on the connected sum. 
\par
Before embarking on the proof of Theorem 0.3, we would like to extend the scope of the invariant $\tilde{s}$ in the following way. Consider a Riemannian spin manifold $(X,g)$ with a Kreck-Stolz product structure (Definition 0.2), that is, suppose that $(X,g)$ is isometric to $(M,g_M) \times (N,g_N)$ (with $(M,g_M)$, $(N,g_N)$ as before), via a spin structure preserving diffeomorphism $X \to M \times N.$ We would like to define $\tilde{s}$ for the manifold $(X,g)$ by declaring $\tilde{s}(X,g):=\tilde{s}(M \times N,g_M+g_N).$ However, we need to argue that such an extension to the definition of $\tilde{s}$ is well-defined, and in order to this we recall:
\proclaim Lemma 2.13. ([EH; page 3075]) If $(X,g)$ is a closed Riemannian manifold, then $X$ decomposes as a Riemannian product of indecomposable factors, and this decomposition is unique in the sense that the corresponding foliations of $X$ are uniquely determined.
\par
In our case, we want to suppose that $(X,g)$ is isometric to the Riemannian product $(M \times N,g_M+g_N).$ In order for our claimed $\tilde{s}$ extension to be well-defined, we need to show that if $(M \times N,g_M+g_N) \cong (M' \times N',g_{M'}+g_{N'})$ then $(M,g_M) \cong (M',g_{M'})$ and $(N,g_N) \cong (N',g_{N'}).$ Now it could be that one or both of $(M,g_M)$ and $(N,g_N)$ are themselves decomposable as a Riemannian product, so we cannot make the desired conclusion directly from Lemma 2.13. However, we will now demonstrate that the topological conditions imposed on $M$ and $N$ in fact provide enough extra structure for us to make this claim.
\proclaim Lemma 2.14. With $M,N$ as before, decompose the manifolds as smooth products $M_1 \times \cdots M_p$ and $N_1 \times \cdots \times N_q$, where the various factors cannot be further decomposed as smooth products. Then no factor of $M$ is homeomorphic to any factor of $N$.
\par
\noindent{\bf Proof.} Let us consider real Pontrjagin classes, since all such classes for $M$ vanish by assumption. Suppose that $N$ splits as a smooth (but not necessarily Riemannian) product $N=N_1 \times \cdots N_q,$ and that one of these factors, $N_1$ say, is also a factor of $M$ up to homeomorphism, i.e. $M\cong N_1 \times K,$ for some $K$. Note that we can work here with homeomorphisms, since real Pontrjagin classes are homeomorphism invariants of smooth manifolds. Now $p_i(M)=\sum_{j+k=i}\pi^*_1p_j(N_1) \cup \pi^*_2 p_k(K),$ where $\pi_1,\pi_2$ indicate the projection maps onto the first, respectively second factors of $N_1 \times K.$ Since $\hat{A}(N) \neq 0$ by assumption and the $\hat{A}$-genus is multiplicative for products, it follows that $\hat{A}(N_1) \neq 0$ also. In particular this means that $N_1$ has some non-vanishing real Pontrjagin classes. Suppose that $p_j(N_1) \neq 0$ for some $j$. We then have $$p_j(M)=\pi^*_1p_j(N_1)+\sum_{r+s=j, r<j}\pi^*_1p_r(N_1) \cup \pi^*_2 p_s(K).$$ Notice that the (non-zero) term $\pi^*_1p_j(N_1)$ belongs to a different summand of $$H^{4j}(M)\cong \bigoplus_{\gamma+\delta=4j}H^\gamma(N_1) \otimes H^\delta(K)$$ than the other terms in the above expression, which have $r<j.$ Hence the other terms cannot cancel out the contribution to $p_j(M)$ coming from $p_j(N_1),$ and we conclude that $p_j(M) \neq 0.$ As this is a contradiction, we deduce that there are no common factors up to homeomorphism between the decompositions of $M$ and $N,$ as claimed. \hfill\sq

\proclaim Corollary 2.15. With $(M,g_M),$ $(N,g_N)$ as before, if $(M \times N,g_M+g_N) \cong (M' \times N',g_{M'}+g_{N'})$ then $(M,g_M) \cong (M',g_{M'})$ and $(N,g_N) \cong (N',g_{N'}).$
\par
\noindent{\bf Proof.} Decompose $(M,g_M)$ and $(N,g_N)$ into Riemannian products with indecomposable factors, and similarly decompose $(M',g_{M'})$ and $(N',g_{N'})$. The image of the isometry $(M \times N,g_M+g_N) \to (M' \times N',g_{M'}+g_{N'})$ provides us with a second isometric splitting of $(M' \times N',g_{M'}+g_{N'})$ as a Riemannian product. By Lemma 2.13 these two splittings must coincide, so our isometry splits as a product of isometries from the factors of $(M \times N,g_M+g_N)$ to the factors of $(M' \times N',g_{M'}+g_{N'}).$ By Lemma 2.14, $M$ and $N$ respectively $M'$ and $N'$ have no common factors, hence we can assemble these factorwise maps into isometries $(M,g_M) \cong (M',g_{M'})$ and $(N,g_N) \cong (N',g_{N'})$ as claimed. \hfill\sq
\medskip
We also note that if the isometry $(M \times N,g_M+g_N) \cong (M' \times N',g_{M'}+g_{N'})$ is orientation preserving, then the isometries $(M,g_M) \cong (M',g_{M'})$ and $(N,g_N) \cong (N',g_{N'})$ are either both orientation preserving or both orientation reversing.

From Corollary 2.15 we immediately deduce:
\proclaim Corollary 2.16. If $(X,g)$ has a Kreck-Stolz product structure with respect to a product $(M \times N,g_M+g_N),$ then the Kreck-Stolz product structure is unique up to isometries of the individual product factors which are either both orientation preserving or both orientation reversing.
\par
Observe that $\tilde{s}$ in Definition 2.8(a) is clearly invariant under orientation preserving isometries of the factors $M$ and $N$, and using the arguments underpinning Lemma 2.11 (see [TW] pages 45-46) also invariant under orientation reversing isometries of both factors. We can therefore now complete the definition of $\tilde{s}$:
\proclaim Definition 2.8(b). Given a closed Riemannian spin manifold $(X,g)$ with a Kreck-Stolz product structure involving an oriented isometry to a product $(M \times N,g_M+g_N)$ as in Definition 2.8(a), we set $\tilde{s}(X,g)=\tilde{s}(M \times N,g_M+g_N),$ where the latter quantity is that defined in Definition 2.8(a).
\par
\noindent{\it Remark:} It is easily observed that Lemmas 2.10 and 2.11 can be immediately extended to incorporate the more general $\tilde{s}$ definition in 2.8(b).
\medskip

\noindent{\bf Proof of Theorem 0.3.} Consider a path of positive scalar curvature metrics $g_{M \times N}(t)$ on $M \times N$ for $t \in [0,1]$ say, where $g_{M \times N}(0)$ and $g_{M \times N}(1)$ are both product metrics with respect to the smooth product structure on $M \times N$. (Note that there is no need to assume that $g_{M\times N}(t)$ is a product metric for any $t \neq 0,1.$) We first establish that $\tilde{s}(M \times N,g_{M \times N}(0))=\tilde{s}(M \times N,g_{M \times N}(1)).$ 

It follows from a well-known observation about paths of positive scalar curvature metrics (see for example [Wr1; Lemma 6.3]) that $g(t)$ can be adjusted to give a metric $g_{M \times N \times I}$ on $M \times N \times I$ for some interval $I$, which has positive scalar curvature globally, agrees with the metrics $g_{M \times N}(0)$ respectively $g_{M \times N}(1)$ when restricted to the two boundary components, and moreover is a product with respect to the $t$ parameter near these boundary components. Thus taking $W=M\times I$, we see that by Lemma 2.9 we have $$\tilde{s}(M \times N \sqcup (-M) \times N,g_{M \times N}(0) \sqcup g_{M \times N}(1))=\hat{A}(N)t(M \times I).$$ By Lemmas 2.10 and 2.11 the left-hand side of this expression is equal to $$\tilde{s}(M \times N,g_{M \times N}(0))-\tilde{s}(M \times N,g_{M \times N}(1)).$$ We claim that $t(M \times I)=0.$ Now the $p_i(M \times I)$ vanish as the Pontrjagin classes of $I$ and (by assumption) the Pontrjagin classes of $M$ both vanish. Thus the $\langle (\hat{A}+a_{k-l}L)(\{j^{-1}p_i(M \times I)\}),[M\times I,\partial (M \times I)]\rangle$ term in $t(M \times I)$ must also be zero. It remains to show that the signature $\sigma(M \times I)=0$, but this follows since $M \times I \simeq M$ and so $H^{4k}(M \times I)=H^{4k}(M^{4k-1})=0.$ Thus we have shown that $\tilde{s}$ is an invariant of product metrics on $M \times N$ belonging to the same path-component of positive scalar curvature metrics.

More generally suppose that both $(X,g)$ and $(X,g')$ both have a Kreck-Stolz product structure involving the same smooth product $M \times N$ and the same diffeomorphism $\phi.$ Let $g(t),$ $t \in [0,1]$ be a path of positive scalar curvature metrics on $X$ with $g(0)=g$ and $g(1)=g'.$ The push-forward metrics $\phi_*(g(t))$ give a path of positive scalar curvature metrics on $M \times N$ beginning with a product metric $g_M+g_N$ and ending with a product metric $g'_M+g'_N.$ According to Definition 2.8(b) we have $\tilde{s}(X,g)=\tilde{s}(M \times N, \phi_*g),$ and by the above paragraph we have $\tilde{s}(M \times N, \phi_*g)=\tilde{s}(M \times N, \phi_*g').$ By Definition 2.8(b) this last term is equal to $\tilde{s}(X,g'),$ and so we deduce that $\tilde{s}(X,g)=\tilde{s}(X,g').$

The assertion 0.3(i), that for $g,g'$ as in the preceding paragraph (though not necessarily in the same path-component of positive scalar curvature metrics), the relative index is given by $i(g,g')=\tilde{s}(X,g)-\tilde{s}(X,g'),$ now follows from equation (\dag) after 2.8(a) in conjunction with the above arguments. In detail, we have $$\tilde{s}(X,g)-\tilde{s}(X,g')=\tilde{s}(M \times N,\phi_\ast(g))-\tilde{s}(M \times N,\phi_\ast(g')).$$ If $g(t)$ is {\it any} smooth path of metrics on $X$ (i.e. with no condition on the scalar curvature) satisfying $g(0)=g,$ $g(1)=g',$ then for the path $\phi_\ast(g(t))$ on $M \times N$ we have $$\eqalign{i(\phi_\ast(g),\phi_\ast(g'))&=\hbox{ind}D^+(M \times N \times I,\phi_\ast(g(t))+dt^2) \cr &=\tilde{s}(M \times N,\phi_\ast(g))-\tilde{s}(M \times N,\phi_\ast(g')),\cr}$$ where the second equality follows from equation (\dag) together with Lemmas 2.10 and 2.11. We also have $$i(g,g')=\hbox{ind}D^+(X \times I,g(t)+dt^2).$$ Now $\phi \times \hbox{id}_I:(X \times I,g(t)+dt^2) \to (M \times N \times I,\phi_\ast(g(t))+dt^2)$ is an orientation preserving isometry. As spin structures are uniquely determined in our circumstances by the orientation, we see that $\phi \times \hbox{id}_I$ is also spin structure preserving. As the index is invariant under spin structure preserving isometries, we deduce that $$i(g,g')=i(\phi_\ast(g),\phi_\ast(g')),$$ and therefore $i(g,g')=\tilde{s}(X,g)-\tilde{s}(X,g')$ as claimed.

To establish assertion 0.3(ii), that the space of positive scalar curvature metrics on $X$ has infinitely many path-components distinguished by $\tilde{s}$, we begin by noting that the argument here is analogous to that of [KS; 2.15]. By [Ca] there is a positive scalar curvature metric $g$ on $S^{4(k-l)-1}$ which is extendable to a positive scalar curvature metric on a certain parallelisable bounding manifold (constructed by plumbing disc bundles). This bounding manifold has non-zero signature and vanishing Pontrjagin classes. It follows from Lemma 2.9 that $\tilde{s}(S^{4(k-l)-1} \times N,g+g_N)$ is a (non-zero) multiple of the (non-zero) signature of the bounding manifold. Consider the manifold $$((M \sharp S_1^{4(k-l)-1}\sharp \cdots \sharp S_p^{4(k-l)-1})\times N,(g_M \sharp g \sharp \cdots \sharp g)+g_N).$$ This is isometric to $(M \times N;g_p+g_N)$ for some positive scalar curvature metric $g_p$ on $M$. Applying Lemma 2.9 to this latter manifold, or Lemma 2.12 to the former, we obtain a different $\tilde{s}$-value for each $p \in \hbox{\bb N}.$ Hence the result in this case.

More generally, consider $(X,g)$ orientation preserving isometric to $(M \times N,g_M+g_N).$ For any $p \in \hbox{\bb N}$ we have a diffeomorphism $M \times N \cong (M \sharp S_1^{4(k-l)-1}\sharp \cdots \sharp S_p^{4(k-l)-1})\times N.$ Composing this isometry and diffeomorphism, then pulling-back the metric $(g_M \sharp g \sharp \cdots \sharp g)+g_N$ to $X$ via this composition, gives a Riemannian manifold $(X,h_p)$ with a Kreck-Stolz product structure. By Definition 2.8(b) we have that $\tilde{s}(X,h_p)=\tilde{s}(M \times N,g_p+g_N),$ which completes the argument in the general case.

Finally, we turn our attention to the moduli space of positive scalar curvature metrics on $X$. First note that the group of diffeomorphisms of $X$ acts on the set of metrics with a Kreck-Stolz product structure, so if $g$ is such a metric, then every representative of its moduli space class $[g]$ has a Kreck-Stolz product structure. Moreover this action is such that all Kreck-Stolz structures belonging to a given orbit involve the same Riemannian product $(M \times N,g_M+g_N).$ If metrics $g$ and $h$ on $X$ differ by an orientation preserving diffeomorphism, it is an immediate consequence of Definition 2.8(b) that $\tilde{s}(X,g)=\tilde{s}(X,h).$ However a diffeomorphism $X \to X$ could reverse orientation, and therefore fail to preserve spin structures. In this case, though, we know from Lemma 2.11 (and the remark after Definition 2.8(b)) that the sign of $\tilde{s}$ changes. Thus we conclude that $|\tilde{s}|$ is invariant under the pull-back action of $\hbox{Diff}(X)$, establishing 0.3(iii). This also gives a well-defined function on the moduli space of Kreck-Stolz structures on $X$, ${\cal K}(X)/\hbox{Diff}(X) \to \hbox{\bb Q},$ since by Corollary 2.16 such a structure is essentially determined by the metric.
\hfill\sq 

\bigskip\bigskip


\centerline{\bf REFERENCES}
\bigskip\bigskip
\item{[APS]} M. F. Atiyah, V.K. Patodi, I. M. Singer, {\it Spectral asymmetry and Riemannian geometry I}, Math. Proc. Camb. Phil. Soc. {\bf 77} (1975), 43-69.
\medskip
\item{[Be]} A.L. Besse, {\it Einstein Manifolds}, Springer-Verlag, Berlin (2002).
\medskip
\item{[BERW]} B. Botvinnik, J. Ebert, O. Randal-Williams, {\it Infinite loop spaces and positive scalar curvature}, arxiv:1411.7408.
\medskip
\item{[BH]} I. Belegradek and J. Hu, {\it Connectedness properties of the space of complete non-negatively curved planes}, Math. Ann. {\bf 362} (2015), 1273-1286. Erratum: Math. Ann. {\bf 364} (2016), 711-712. 
\medskip
\item{[BHSW]} B. Botvinnik, B. Hanke, T. Schick, M. Walsh, {\it Homotopy groups of the moduli space of metrics of positive scalar curvature}, Geom. Topol. {\bf 14} (2010), 2047-2076.
\medskip
\item{[CM]} F. Cod\'a Marques, {\it Deforming three-manifolds with positive scalar curvature}, Ann. of Math. (2) {\bf 176} (2012), no. 2, 815-863.
\medskip 
\item{[CS]} D. Crowley, T. Schick, {\it The Gromoll filtration, $KO$-characteristic classes and metrics of positive scalar curvature}, Geom. Topol. {\bf 17} (2013), 1773-1790.
\medskip
\item{[DKT]} A. Dessai, S. Klaus, W. Tuschmann, {\it Nonconnected moduli spaces of nonnegative sectional curvature metrics on simply-connected manifolds}, Bull. London Math. Soc. {\bf 50} (2018), 96-107.
\medskip
\item{[EH]} J.-H. Eschenburg, E. Heintze, {\it Unique decomposition of Riemannian manifolds}, Proc. Amer. Math. Soc. {\bf 126}, no. 10, (1998), 3075-3078.
\medskip
\item{[FO1]} F. T. Farrell, P. Ontaneda, {\it The Teichm\"uller space of pinched negatively curved metrics on a hyperbolic manifold is not contractible}, Ann. of Math. (2) {\bf 170} (2009), no. 1, 45-65.
\medskip
\item{[FO2]} F. T. Farrell, P. Ontaneda, {\it The moduli space of negatively curved metrics of a hyperbolic manifold}, J. Topol. {\bf 3} (2010), no. 3, 561-577.
\medskip
\item{[FO3]} F. T. Farrell, P. Ontaneda, {\it On the topology of the space of negatively curved metrics}, 
J. Diff. Geom. {\bf 86} (2010), no. 2, 273-301.
\medskip
\item{[GL]} M. Gromov, H.B. Lawson, {\it The classification of manifolds of positive scalar curvature}, Ann. Math. {\bf 111} (1980), 423-434.
\medskip
\item{[HSS]} B. Hanke, T. Schick, W. Steimle, {\it The space of metrics of positive scalar curvature}, Publ. Math. Inst. Hautes \'Etudes Sci. {\bf 120} (2014), 335-367.
\medskip
\item{[KPT]} V. Kapovitch, A. Petrunin, W. Tuschmann, {\it Non-negative pinching, moduli spaces and bundles with
infinitely many souls}, J. Diff. Geom. {\bf 71} (2005) no. 3, 365-383. 
\medskip
\item{[KS]} M. Kreck, S. Stolz, {\it Nonconnected moduli spaces of positive sectional curvature metrics}, J. Am. Math. Soc. {\bf 6} (1993), 825-850.
\medskip
\item{[LM]} H.B. Lawson, M.-L. Michelsohn, {\it Spin Geometry}, Princeton Math. Series {\bf 38}, Princeton University Press, (1989).
\medskip
\item{[Mo]} S. Morita, {\it Geometry of Differential Forms}, Translations of Mathematical Monographs vol. {\bf 201}, American Mathematical Society (2001).
\medskip
\item{[MT]} I. Madsen, J. Tornehave, {\it From Calculus to Cohomology}, Cambridge University Press (1997).
\medskip
\item{[P]} P. Petersen, {\it Riemannian Geometry}, Graduate Texts in Mathematics {\bf 171}, Springer-Verlag, (1998).
\medskip
\item{[PS]} P. Piazza, T. Schick, {\it Groups with torsion, bordism and rho-invariants}, Pacific J. Math. {\bf 232} (2007), no. 2, 355-378.
\medskip
\item{[R]} J. Rosenberg, {\it Manifolds of positive scalar curvature: a progress report}, Surveys in Differential Geometry vol. XI, International Press (2007), 259-294.
\medskip
\item{[RS]} J. Rosenberg, S. Stolz, {\it Metrics of positive scalar curvature and connections with surgery}, Surveys in Surgery Theory vol. 2, Annals of Mathematics Studies {\bf 149}, Princeton University Press, (2001).
\medskip
\item{[S1]} T. Schick, {\it A counterexample to the unstable Gromov-Lawson-Rosenberg conjecture}, Topology {\bf 37} (1998), 1165-1168.
\medskip
\item{[S2]} T. Schick, {\it The topology of scalar curvature}, arXiv:1405.4220v2.
\medskip
\item{[TW]} W. Tuschmann, D. J. Wraith, {\it Moduli spaces of Riemannian metrics}, Oberwolfach Seminars {\bf 46}, Birkh\"auser, Springer Basel (2015).
\medskip
\item{[TWi]} W. Tuschmann, M. Wiemeler, {\it On the topology of moduli spaces of non-negatively curved Riemannian metrics}, arXiv:1712.07052.
\medskip
\item{[Wa1]} M. Walsh, {\it Cobordism invariance of the homotopy type of the space of positive scalar curvature metrics}, Proc. Amer. Math. Soc. {\bf 141} (2013), no. 7, 2475-2484.
\medskip
\item{[Wa2]} M. Walsh, {\it H-spaces, loop spaces and the space of positive scalar curvature metrics on the sphere}, Geom. Topol. {\bf 18} (2014), no. 4, 2189-2243.
\medskip
\item{[Wr1]} D. J. Wraith, {\it On the moduli space of positive Ricci curvature metrics on homotopy spheres}, Geom. Topol. {\bf 15} (2011), 1983-2015.
\medskip
\item{[Wr2]} D. J. Wraith, {\it Non-negative versus positive scalar curvature}, arXiv:1607.00657.
\bigskip\medskip
\noindent {\it Department of Mathematics and Statistics, National University of Ireland Maynooth,}

\noindent {\it Maynooth, County Kildare, Ireland. Email: david.wraith@mu.ie.}

\end